\documentclass[11pt]{article}
\usepackage{amssymb}

\begin{document}

\newcommand{\ls}[1]
   {\dimen0=\fontdimen6\the\font \lineskip=#1\dimen0
\advance\lineskip.5\fontdimen5\the\font \advance\lineskip-\dimen0
\lineskiplimit=.9\lineskip \baselineskip=\lineskip
\advance\baselineskip\dimen0 \normallineskip\lineskip
\normallineskiplimit\lineskiplimit \normalbaselineskip\baselineskip
\ignorespaces }

\newtheorem{theorem}{Theorem}
\newtheorem{problem}{Problem}
\newtheorem{definition}{Definition}
\newtheorem{lemma}{Lemma}
\newtheorem{proposition}{Proposition}
\newtheorem{corollary}{Corollary}
\newtheorem{example}{Example}
\newtheorem{conjecture}{Conjecture}
\newtheorem{algorithm}{Algorithm}
\newtheorem{exercise}{Exercise}
\newtheorem{remarkk}{Remark}

\newcommand{\be}{\begin{equation}}
\newcommand{\ee}{\end{equation}}
\newcommand{\bea}{\begin{eqnarray}}
\newcommand{\eea}{\end{eqnarray}}
\newcommand{\beq}[1]{\begin{equation}\label{#1}}
\newcommand{\eeq}{\end{equation}}
\newcommand{\beqn}[1]{\begin{eqnarray}\label{#1}}
\newcommand{\eeqn}{\end{eqnarray}}
\newcommand{\beaa}{\begin{eqnarray*}}
\newcommand{\eeaa}{\end{eqnarray*}}
\newcommand{\req}[1]{(\ref{#1})}

\newcommand{\lip}{\langle}
\newcommand{\rip}{\rangle}
\newcommand{\uu}{\underline}
\newcommand{\oo}{\overline}
\newcommand{\La}{\Lambda}
\newcommand{\la}{\lambda}
\newcommand{\eps}{\varepsilon}
\newcommand{\om}{\omega}
\newcommand{\Om}{\Omega}

\renewcommand{\lll}{{\Bigl(}}
\newcommand{\rrr}{{\Bigr)}}
\newcommand{\qqq}{{\Bigl\|}}

\newcommand{\dint}{\displaystyle\int}
\newcommand{\dsum}{\displaystyle\sum}
\newcommand{\dfr}{\displaystyle\frac}
\newcommand{\bige}{\mbox{\Large\it e}}
\newcommand{\integers}{{\mathbb Z}}
\newcommand{\rationals}{{\Bbb Q}}
\newcommand{\reals}{{\rm I\!R}}
\newcommand{\realsd}{\reals^d}
\newcommand{\realsn}{\reals^n}
\newcommand{\NN}{{\rm I\!N}}
\newcommand{\DD}{{\rm I\!D}}
\newcommand{\degree}{{\scriptscriptstyle \circ }}
\newcommand{\dfn}{\stackrel{\triangle}{=}}
\def\complex{\mathop{\raise .45ex\hbox{${\bf\scriptstyle{|}}$}
     \kern -0.40em {\rm \textstyle{C}}}\nolimits}
\def\hilbert{\mathop{\raise .21ex\hbox{$\bigcirc$}}\kern -1.005em {\rm\textstyle{H}}} 
\newcommand{\RAISE}{{\:\raisebox{.6ex}{$\scriptstyle{>}$}\raisebox{-.3ex}
           {$\scriptstyle{\!\!\!\!\!<}\:$}}} 

\newcommand{\hh}{{\:\raisebox{1.8ex}{$\scriptstyle{\degree}$}\raisebox{.0ex}
           {$\textstyle{\!\!\!\! H}$}}}

\newcommand{\OO}{\won}
\newcommand{\calA}{{\cal A}}
\newcommand{\calB}{{\cal B}}
\newcommand{\calC}{{\cal C}}
\newcommand{\calD}{{\cal D}}
\newcommand{\calE}{{\cal E}}
\newcommand{\calF}{{\cal F}}
\newcommand{\calG}{{\cal G}}
\newcommand{\calH}{{\cal H}}
\newcommand{\calK}{{\cal K}}
\newcommand{\calL}{{\cal L}}
\newcommand{\calM}{{\cal M}}
\newcommand{\calO}{{\cal O}}
\newcommand{\calP}{{\cal P}}
\newcommand{\calX}{{\cal X}}
\newcommand{\calXX}{{\cal X\mbox{\raisebox{.3ex}{$\!\!\!\!\!-$}}}}
\newcommand{\calXXX}{{\cal X\!\!\!\!\!-}}
\newcommand{\gi}{{\raisebox{.0ex}{$\scriptscriptstyle{\cal X}$}
\raisebox{.1ex} {$\scriptstyle{\!\!\!\!-}\:$}}}
\newcommand{\intsim}{\int_0^1\!\!\!\!\!\!\!\!\!\sim}
\newcommand{\intsimt}{\int_0^t\!\!\!\!\!\!\!\!\!\sim}
\newcommand{\pp}{{\partial}}
\newcommand{\al}{{\alpha}}
\newcommand{\sB}{{\cal B}}
\newcommand{\sL}{{\cal L}}
\newcommand{\sF}{{\cal F}}
\newcommand{\sE}{{\cal E}}
\newcommand{\sX}{{\cal X}}
\newcommand{\R}{{\rm I\!R}}
\newcommand{\vp}{\varphi}
\newcommand{\N}{{\rm I\!N}}
\def\ooo{\lip}
\def\ccc{\rip}
\newcommand{\ot}{\hat\otimes}
\newcommand{\rP}{{\Bbb P}}
\newcommand{\bfcdot}{{\mbox{\boldmath$\cdot$}}}

\renewcommand{\varrho}{{\ell}}
\newcommand{\dett}{{\textstyle{\det_2}}}
\newcommand{\sign}{{\mbox{\rm sign}}}
\newcommand{\TE}{{\rm TE}}
\newcommand{\TA}{{\rm TA}}
\newcommand{\E}{{\rm E\,}}
\newcommand{\won}{{\mbox{\bf 1}}}
\newcommand{\Lebn}{{\rm Leb}_n}
\newcommand{\Prob}{{\rm Prob\,}}
\newcommand{\mod}{{\rm mod\,}}
\newcommand{\sinc}{{\rm sinc\,}}
\newcommand{\ctg}{{\rm ctg\,}}
\newcommand{\loc}{{\rm loc}}
\newcommand{\trace}{{\,\,\rm trace\,\,}}
\newcommand{\Dom}{{\rm Dom}}
\newcommand{\ifff}{\mbox{\ if and only if\ }}
\newcommand{\proof}{\noindent {\bf Proof:\ }}
\newcommand{\remark}{\noindent {\bf Remark:\ }}
\newcommand{\remarks}{\noindent {\bf Remarks:\ }}
\newcommand{\note}{\noindent {\bf Note:\ }}

\newcommand{\boldx}{{\bf x}}
\newcommand{\boldX}{{\bf X}}
\newcommand{\boldy}{{\bf y}}
\newcommand{\boldR}{{\bf R}}
\newcommand{\uux}{\uu{x}}
\newcommand{\uuY}{\uu{Y}}

\newcommand{\limn}{\lim_{n \rightarrow \infty}}
\newcommand{\limN}{\lim_{N \rightarrow \infty}}
\newcommand{\limr}{\lim_{r \rightarrow \infty}}
\newcommand{\limd}{\lim_{\delta \rightarrow \infty}}
\newcommand{\limM}{\lim_{M \rightarrow \infty}}
\newcommand{\limsupn}{\limsup_{n \rightarrow \infty}}

\newcommand{\ra}{ \rightarrow }

\newcommand{\ARROW}[1]
  {\begin{array}[t]{c}  \longrightarrow \\[-0.2cm] \textstyle{#1} \end{array} }

\newcommand{\AR}
 {\begin{array}[t]{c}
  \longrightarrow \\[-0.3cm]
  \scriptstyle {n\rightarrow \infty}
  \end{array}}

\newcommand{\pile}[2]
  {\left( \begin{array}{c}  {#1}\\[-0.2cm] {#2} \end{array} \right) }

\newcommand{\floor}[1]{\left\lfloor #1 \right\rfloor}

\newcommand{\mmbox}[1]{\mbox{\scriptsize{#1}}}

\newcommand{\ffrac}[2]
  {\left( \frac{#1}{#2} \right)}

\newcommand{\one}{\frac{1}{n}\:}
\newcommand{\half}{\frac{1}{2}\:}

\def\le{\leq}
\def\ge{\geq}
\def\lt{<}
\def\gt{>}

\def\squarebox#1{\hbox to #1{\hfill\vbox to #1{\vfill}}}
\newcommand{\qed}{\hspace*{\fill}
           \vbox{\hrule\hbox{\vrule\squarebox{.667em}\vrule}\hrule}\bigskip}

\title{Some measure-preserving  point transformations on the Wiener
space and their ergodicity}

\author{A.S. \"Ust\"unel and M. Zakai}
\date{}
\maketitle

\noindent
{\bf Abstract:} 
{\small{
Suppose that $T$ is a map of the Wiener space into itself, of the
following type: $T=I+u$ where $u$ takes its values in the
Cameron-Martin space $H$. Assume also that $u$ is a finite sum of
$H$-valued multiple Ito-Wiener integrals. In this work we prove that
if $T$ preserves the Wiener measure, then necessarily $u$ is in the
first Wiener chaos and the transformation corresponding  to it is  a
rotation in the sense of \cite{USZ}. 
 Afterwards the ergodicity and mixing  of rotations which
 are second  quantizations of the unitary operators on  the Cameron-Martin
    space, are  characterized. Finally, the ergodicity of the
transformation $dY_t = \gamma (t) d W_t$, $0\le t \le 1$
where $W$ is $n$-dimensional Wiener and $\gamma$ is non random is
characterized.
}}  

\section{Introduction}

Let $\mu $ be the  standard  Gaussian measure on $\reals^n$,  i.e.
$$
\mu \Bigl\{x:x_i\leq  a_i, i=1,2, \ldots,n\Bigr\} = \prod_{i=1}^n \Phi(a_i)
$$ 
where 
\beq{1}
\Phi (a) = (2\pi)^{-\half} \int_{-\infty}^a e ^{-\frac{\eta^2}{2}} 
d \eta\,.
\eeq
Then
\begin{itemize}
\item[(a)]The linear point-transformations $T$ on $\reals^n$ which leave
this measure invariant induce  unitary transformations on $L^2(\mu,\R^n)$, 
which are defined as $Of(x)=f\circ T(x)$.
\item[(b)]  There are many non-linear transformations on $\reals^n$ 
which leave the  measure $\mu$  invariant, too many to characterize without
any further restriction. 
\item[(c)]The   transformation   $T$  is not
ergodic: in fact let $f$ be defined as $f(x)=|x|_{\R^n}$, then
$$
f=f\circ T
$$
and evidently, $f$ is a non-constant function.
\end{itemize}
\noindent
The infinite dimensional extension of this problem leads directly to
the formulation of the problem for Wiener processes.  Indeed, let
$(w_t, t \in [0,1])$ denote the standard Wiener process and let
$(e_i ,\,i\in \NN)$
be  a complete orthonormal basis  in the Cameron-Martin space $H$. Denote by
$(e'_i,i\in \NN)$ the image of this basis in $L_2([0,1])$. Define 
\beq{2}
\delta e_i = \int_0^1 e'_i(s) dw(s)\,,
\eeq
then
$(\delta e_i,\,i\in \NN)$ are i.i.d. $N(0,1)$-random variables and 
\cite{ITN}, \cite{BOOK}
\beq{3}
w_t = \sum_1^\infty \delta e_i  e_i(t)
 = \sum_1^\infty \delta e_i \int_0^t e'_i(s) ds
\eeq
in the sense that
$$
\sup_{t\in[0,1]} \left| w_t - \sum_1^N \delta e_i \int_0^t e'_i (s) ds
\right|
\begin{array}{c} 
{\rm a.s.}\\[-.2cm] \longrightarrow \\[-.2cm] N\to\infty
\end{array} 
0\,.
$$
It follows that for any transformation, linear or non-linear, invertible or
non-invertible, from $\{\delta e_i, i=1,2,\cdots\}$ to another sequence,
say $(\eta_i,\, i\in \NN)$, of i.i.d. $N(0,1)$-random variables,
defined by 
$$
\sum_{i=1}^\infty \delta e_i e_i\to 
\sum_{i=1}^\infty \eta_i  e_i 
$$
will be a measure invariant transformation of the Wiener space.

A class of transformations which plays an important role in many
applications is the shift transformation
\beq{4}
(Tw)_t = w_t + \int_0^t u_s (w) ds, 
0 \le t \le 1,
\eeq
where
\beq{5}
\int_0^1|u_s|^2 ds < \infty   \quad {\rm a.s.}
\eeq
It is natural to ask for a characterization of the shifts $u$ for which
$T$ is measure-invariant, i.e. $Tw$ is also a Wiener process on
$C_0([0,1])$.  In the next section we consider the transformations induced
by a finite sum of multiple Wiener-Ito integrals taking values in the
Cameron-Martin space and characterize those
shifts which induce an invariant measure. We prove in particular their 
non-ergodicity.
In section~3 we study the measure preserving transformations which are 
defined via the second quantization of deterministic  unitary operators on the
Cameron-Martin space which cover also the special kind of shifts
presented in the second section. 
In particular a necessary and
sufficient condition for their ergodicity and mixing 
is proved.
Section~4 deals with the special case where 
$$
dY_t=\gamma(t) d W_t, \quad 0 \le t \le 1\,, \quad Y_0=0\,, 
$$
where $W$ is a standard
$n$-dimensional Wiener process and $\gamma(t)$ is not random and takes
values in the group of unitary matrices.  The ergodicity of this
transformation is characterized.

The characterization of ergodicity and mixing for real valued Gaussian
processes is due to Maruyama (cf.\cite{MAR}).  The results presented
in Theorems \ref{ergos-thm} and  Theorem \ref{mixing-thm}  are infinite
dimensional extensions of the results of 
Maruyama  and can be derived by starting from Maruyama's results
(bypassing Lemma~\ref{conty-lemma}).  We preferred, however, the proof
presented here as it is 
more direct and shorter.  It is based on the following 
characterizations (cf. e.g.\
section~1.7 of \cite{COF}).
Let $T$ be an automorphism (invertible, $T$ and $T^{-1}$ are measurable and
measure preserving) then:
\begin{itemize}
\item[(A)]
$T$ is ergodic, if and
only if the only eigenfunctions of the induced unitary transformation 
$O$ associated with $\la=1$ are the constants.
\item[(B)]
$T$ is weak mixing if and only if $O$ has no eigenfunctions other than
constants.
\item[(C)]
Let $L^2_0(\mu)$ denote the class of real
valued square integrable, zero mean Wiener functionals.
Set $a_n(f) = E[(O^n f)\cdot f]$.  Then $T$ is mixing if and only if 
$a_n(f) \to 0$ as
$n\to \infty$ for all $f$ in $L_0^2(\mu)$.
\end{itemize}

\noindent
{\bf Remarks}

\noindent
(a)~The results presented here are valid
for arbitrary abstract Wiener spaces although the study here is in  the
setup of the classical Wiener space.

\noindent
(b)~The ergodicity problem considered in this paper deals with invertible
transformations.  The invertibility however, is not necessary for
ergodicity.  Indeed, let $w_\bfcdot$ be as in equation \req{3}, set
$$
(Tw)_t = \sum_1^\infty \delta e_{i+1} e_i(t)
$$
then it is easily verified that $T$ is measure preserving and strong
mixing.

\noindent
(c)~After this paper was written we learned of the paper \cite{11} by
Wiener and Akutowicz which characterizes the mixing properties of
transformation discussed in section~3.

\section{Shifts induced by multiple Wiener-Ito\newline integrals}
In the sequel we denote by $(C_0([0,1]),H,\mu)$ the classical Wiener
space, where $H$ denotes the Cameron-Martin space which consists of
absolutely continuous functions on $[0,1]$ with square integrable
derivatives and $\mu$ is the Wiener measure. Recall that one can
define a Sobolev derivative on this 
space respecting the $\mu$-equivalence classes (cf. e.g. \cite{ASU}),
whose adjoint, denoted by $\delta$, called divergence operator, which
coincides with the Ito integral of the Lebesgue density of $H$-valued
functional if the latter is adapted to the filtration of the Wiener
process. 
Let $\{w_t, t\in [0,1]\}$ be the standard Wiener process on
$C_0([0,1])$. Assume that  
$k_{n+1}\in L^2([0,1]^{n+1})$ is a kernel which is 
symmetric in its  first $n$ variables.  Let
$I_n (k_{n+1} (s_1, \cdots, s_n, t))$ or just
$I_n (k_{n+1} (\cdot, t))$ denote the $n$-th order multiple Wiener-Ito
integral with respect to  $s_1,\ldots,s_n$   of $k_{n+1}$. For
$t\in [0,1]$, define 
\beq{6}
y_t = (Tw)_t =w_t + \sum_1^N \int_0^t I_n \Bigl(k_{n+1} ( \cdot, \eta)\Bigr)
d \eta
\eeq
for some finite $N$.  Let $\mu$ be the standard Wiener measure  and
denote by  $T^*\mu$  the measures induced on
$C_0([0,1])$ by $w \to Tw$.    
\begin{theorem}
\label{char-thm}
Let $Tw$ be as defined by \req{6}, then $T^*\mu=\mu$ and only if
\begin{itemize}
\item[(a)]~$N=1$
\item[(b)] and $(I+K)$ is a unitary operator on $L^2([0,1])$, i.e.
$$
(I+K)(I+K)^* = (I+K)^* (I+K) = I\,,
$$
where $K$ is defined on $L^2([0,1])$ by 
$$
Kf(t)=\int_0^1 k_2(t,\tau)f(\tau)d\tau\,.
$$
\end{itemize}
\end{theorem}
\remarks
Condition (b) can be restated as:
\begin{itemize} 
\item[(b$_1$)]$ -1$ is not an eigenvalue of
$K$. 
\item[(b2)] 
$$
k_2(s,t) + k_2(t,s) + \int_0^1 k_2 (\theta,s) k_2 (\theta, t) d\theta=0
$$
or equivalently
$$
k_2(s,t) + k_2(t,s) + \int_0^1 k_2 (s,\theta) k_2 (t,\theta) d\theta=0\,.
$$
for any $(s,t)\in [0,1]^2$,  $ds\times dt$ almost surely.
\end{itemize}
\proof
To show necessity, let $h(t)$ be in $L^2([0,1])$.
If $T^*\mu=\mu$ then
$$
\int_0^1 h(s) dy_s = \int_0^1 h(s) dw_s +\sum_{n=1}^N \int_0^1 h(\eta) I_N
\Bigl(k_{n+1} (\cdot, \eta)\Bigr) d \eta
$$
is a zero mean Gaussian random variable.  By a standard convergence argument,
the order of integration can be interchanged and it holds that
\beq{7}
\int_0^1 h(s) dy_s = \int_0^1 h(s) dw_s + \sum_1^N I_n 
\left( \int_0^1 k_{n+1} ( \cdot, \eta)  h (\eta) d\eta\right)\,.
\eeq
The term on the left hand side is Gaussian and for $n\ge 2$,
$I_n(\cdot)$ is non-Gaussian.  Moreover, a result of  McKean
(cf.\ section~8 of \cite{McK}) states that if $f_k(s_1, \cdots s_k)$ are
non-zero elements of $L^2([0,1]^k)$ then for some positive $\alpha$ and
$\beta$ and for $x$ large enough
$$
\exp - \alpha x^{2/N} \le \Prob \left(\left| \sum_1^N
I_k (f_k) \right| > x\right) \le \exp - \beta x^{2/N}\,.
$$
Since  there can be no cancellation between the terms in \req{7},
we must have $N=1$, and \req{7} becomes
\beq{8}
\int_0^1 h(s) dy_s = \int_0^1 h(s) dw_s+ \int_0^1
\left(\int_0^1 k_2(s, \theta) h (\theta) d \theta\right) dw_s\,.
\eeq
The operator $K$ corresponding to the kernel  $k_2$ is
Hilbert-Schmidt on $L^2[0,1]$, hence it has a discrete spectrum.
If $\la=-1$ is an eigenvalue of $K$ and $h$ is a
corresponding eigenfunction then, almost surely, 
$\int_0^1 h(s) dy_s=0$ which contradicts the assumption that $w\to y(w)$
is Wiener, 
this yields  condition (b1).  Furthermore, if $w\to y(w)$ is Wiener then
$$
E\left[\int_0^1 g_1(s) dy_s \int_0^1 g_2(s) dy_s\right]=
\int_0^1 g_1(s) g_2(s) ds\,.
$$
Hence, for any $h,\alpha \in H$, by \req{8}
\beaa
E[(\delta h\circ T)\,(\delta\alpha \circ
T)]&=&(h,\alpha)_H+\Bigl(K^*h,K^*\alpha\Bigr)_H\\ 
&&\,+\Bigl(K^*hh,\alpha\Bigr)_H+\Bigl(K^*h,\alpha\Bigr)_H\\
&=&(h,\alpha)_H
\eeaa
hence
$$
\Bigl(KK^*h+K^*h+Kh,\alpha\Bigr)_H=0
$$
therefore (b) and (b2) follow.

\qed

\begin{corollary}
Under the hypothesis of Theorem \ref{char-thm},  the mapping  $T$ is
almost surely invertible \footnote{This means the existence of a
  measurable map $S:W\to W$ such that $\mu\{T\circ S=S\circ
  T=I_W\}=1$, cf. \cite{BOOK}.} and we have  
$$
|\dett(I_H+K)|\exp\left\{-I_2(k_2)-1/2|\delta
  K|_H^2\right\}=1
$$
and 
$$
|\dett(I_H+K)|=1\,,
$$
where $\dett(I_H+K)$ denotes the modified Carleman-Fredholm
determinant (cf. \cite{DUS}).
\end{corollary}
\proof
The hypothesis implies that $T$ is invertible. Indeed, $w\to
T^{-1}(w)$ is given by  
$$
T^{-1}(w)=w-\int_0^\cdot I_1(\beta(t,\cdot))dt\,,
$$
where $\beta(s,t)$ is the symmetric kernel associated to the
Hilbert-Schmidt  operator $(I_H+K)^{-1}K$ (cf. \cite{DeU}).
By  the change of variables formula, for any continuous and bounded
function $f$ on the Wiener space, we have (\cite{KUS}, \cite{BOOK})
$T^*\mu \sim \mu$ and
$$
E[f\circ T\,|\La|]=E[f]\,,
$$
where  
$$
\La= \dett(I_H+K)\exp\left\{-I_2(k_2)
- \half\int_0^1 \left( \int_0^1 k_2(s,t) dw_s\right)^2 dt\right\}\,.
$$
Since $E[|\La|] = 1$, in order to show that
$|\La|= 1 $ it suffices to show that
$$
-I_2(k_2) - \half \int_0^1 \left(\int_0^1 k_2(s,t) dw_s\right)^2 dt
$$
is independent of $w$.  Now, by Ito's rule
\beaa
\lefteqn{
I_2(k_2) + \half \int_0^1 \left( \int_0^1 k_2(s,t) dw_s \right)^2 dt}\\
&& = I_2(k_2) + I_2 \left(\int_0^1 k_2(s,t) k_2(\theta, t) dt\right)
+ \half \int_0^1 \int_0^1 k_2^2 (s,t) ds\,dt
\eeaa
and  the $I_2$ terms must vanish since from $(b_2)$, we have 
$$
k_2(s,\theta) + k_2(\theta,s)+ \int_0^1 k_2(s,t) k_2(\theta,t) dt = 0
$$
and the proof  follows.
\qed

\begin{corollary}
The class of transformations $T$ satisfying the conditions of the
Theorem \ref{char-thm}
form a subgroup of the group of transformations
$$
Tw=w+\int_0^\bfcdot a_s(w) ds\,,
$$
with $\int_0^1 |a_s(w)|^2 ds < \infty$ a.s. for which $T^*\mu=\mu$.
\end{corollary}

\proof
Setting
\beaa
T_1 w(t) & = & w_t + \int_0^t \int_0^1 k(s,\theta) d w_s\, d\theta\\
T_2 w(t) & = & w_t + \int_0^t \int_0^1 q(s,\theta) d w_s\, d\theta
\eeaa
and assuming that $k$ and $q$ satisfy the conditions of the theorem then
$(T_2T_1)^* \mu=\mu$.  Now,
\beaa
T_2(T_1 w)(t) & = & w_t + \int_0^t\int_0^1 k(s,\theta) dw_s d\theta \\
&& + \int_0^t \left[\int_0^1 q(s,\theta) dw_s + \int_0^1 k (\rho,\theta)
dw_\rho \cdot ds \right] d\theta \\
& = & w_t + \int_0^t \int_0^1 k(s,\theta) dw_s d\theta + \int_0^t\int_0^1
q(s,\theta) dw_s d\theta \\
&& + \int_0^t \left( \int_0^1 
\left( \int_0^1 q(\theta,\eta) k(s,\eta) d\eta\right)
dw_s\right) d\theta
\eeaa
and the result follows since $q$ and $k$ are Hilbert-Schmidt kernels, so is
$\int_0^1q(\cdot, \eta) k(\cdot, \eta) d\eta$.
\qed

Such a transformation is never ergodic as it is proven in the following
\begin{proposition}
\label{non-ergos}
Any  transformation of the Wiener space   satisfying the
conditions of Theorem \ref{char-thm}    is non-ergodic.
\end{proposition}

\proof
Assume that $\la$ is an eigenvalue of $K$, with the corresponding
eigenfunction $h$. Then $I_1(h)$ is an eigenfunction of $O$ with the
eigenvalue $1+\la$. Since $O$ is an isometry, we should have
necessarily $|1+\la|=1$, moreover 
\beaa
O|I_1(h)|&=&|I_1(h)\circ T|\\
   &=&|1+\la||I_1(h)|\\
  &=&|I_1(h)|\,.
\eeaa
Consequently, $ |I_1(h)|$ is a non-trivial invariant  function, hence
$f\to f\circ T$  can not be ergodic.

\qed

\section{Ergodicity of transformations induced by\newline rotations }
Let $w$ denote, as before,  the standard Wiener path  and let $R$ be a 
non-random,  unitary  transformation of the Cameron-Martin space $H$.
Let $(e_i, i\in \NN)$ be a complete, orthonormal basis  of $H$ whose
image in  $L^2([0,2\pi])$ will be denoted by $(e'_i)$.  Set 
\beq{4.1}
Tw=\sum_{i=1}^\infty \delta (R e_i) \,e_i\,.
\eeq
Since $\delta (R e_i)$ are i.i.d. and $N(0,1)$, $Tw$ is also a  Wiener path
(cf.\ \cite{USZ} or \cite{BOOK} for more general cases).

\begin{lemma}
\label{chaos-lemma}
The definition of $Tw$ is independent of the choice of the basis 
$( e_i, i\in \NN)$. Hence we have also
\beq{42}
Tw = \sum_1^\infty \delta e_i \cdot (R^{-1} e_i)\,.
\eeq
Moreover, for any $h\in H$, one has
$$
\exp\left\{\delta h-1/2|h|_H^2\right\}\circ T=
\exp\left\{\delta(Rh)-1/2|h|_H^2\right\}\,.
$$
In particular, if $F\in L^2(\mu)$ has the Wiener chaos representation
as 
$$
F=E[F]+\sum_{i=1}^\infty I_n(f_n)\,,
$$
then
\begin{equation}
\label{chaos-dec}
F\circ T=E[F]+\sum_{i=1}^\infty I_n\left(R^{\otimes n}(f_n)\right)\,,
\end{equation}
where $R^{\otimes n}$ denotes $n$-th tensor power of the operator $R$.
\end{lemma}
\proof
Let $\alpha$ be an element of the continuous dual $C([0,1])^\star$   of 
$C([0,1])$, i.e., a 
bounded Borel measure on $[0,1]$. Let $\tilde{\alpha}$ denote its image under 
the injection $C([0,1])^\star \hookrightarrow H$. Then it is easy to see that 
$\tilde{\alpha}(t)=\int_0^t\alpha([s,1])ds$. We have 
\beaa
<T w,\alpha>&=&\sum_{i=1}^\infty \delta(Re_i)\int_0^1e_i(s)d\alpha(s)\\
   &=&\sum_{i=1}^\infty \delta(Re_i)({\tilde{\alpha}},e_i)_H\\
   &=&\sum_{i=1}^\infty \delta(Re_i)(R{\tilde{\alpha}},Re_i)_H\\
   &=&\delta(R{\tilde{\alpha}})
\eeaa
since $(Re_i,i\in \NN)$ is a complete, orthonormal basis of $H$. By the density 
of $C([0,1])^\star$ in $H$, we obtain that $\delta h\circ
T=\delta(Rh)$ for any $h\in H$, the second claim  is now obvious and
the identity (\ref{chaos-dec}) follows from it.
\qed

\remarks
(i)~since $R$ is unitary, it possesses the spectral representation
$$
R= \int_{0}^{2\pi} e^{i\theta} d_\theta \Pi_\theta
$$
where $(\Pi_\theta,\theta\in [0,2\pi])$ is a resolution of the
identity. 
It follows by standard arguments that
$$
\delta h\circ T^n(w)= \int_{0}^{2\pi} e^{in\theta}
d_\theta\delta(\Pi_\theta h)\,,
$$
where the integral at the right hand side is to be interpreted as a
stochastic integral with respect to the martingale $\theta\to \delta
\Pi_\theta h$ (cf. \cite{BOOK}).  Hence 
$$
E\Bigl[ \delta h_1\circ T^n\,  \delta h_2\Bigr] = \int_{0}^{2\pi}
e^{in\theta} d_\theta( \Pi_\theta h_1, h_2)\,.
$$
(ii) The transformation studied   in the previous section (cf. Theorem
\ref{char-thm})  is a particular case 
of this one defined by \req{4.1}  since $I_H+K$ is a unitary operator
on $H$. 

Before proceeding further let us prove a technical result which will
be useful in the sequel: 
 
\begin{lemma}
\label{conty-lemma}
Let $(\Pi_\theta,\theta \in [0,2\pi])$ be   a  resolution of identity on
$H$. Assume that $\theta\to (\Pi_\theta h,k)_H$ is continuous for any 
$h,k\in H$. Then for any $f,g\in H^{\odot n}$ (i.e. the symmetric
tensor product of order $n$),   
$$
(\theta_1,\ldots,\theta_n)
\to\Bigl((\Pi_{\theta_1}\otimes\ldots\otimes
\Pi_{\theta_n})f,g\Bigr)_{H^{\odot n}} 
$$
is continuous on $[0,2\pi]^n$. Moreover, for any $f\in H^{\odot n}$,
$$
(\theta_1,\ldots,\theta_n) \to d((\Pi_{\theta_1}\otimes\ldots
\Pi_{\theta_n})f,f)_{H^{\odot n}}
$$
is a $\sigma$-additive and  atomless  measure on $[0,2\pi]^n$.
\end{lemma}

\proof
If  $f=a_1\odot\ldots\odot a_n$ and  
$g=h_1\odot\ldots\odot h_n$, with $a_i, h_j\in H$ (we say in this case
that $f$ and $g$ are pure vectors), then, denoting the vector
$(\theta_1,\ldots,\theta_n)$ by $\vec{\theta}$,  
$$
(\Pi_{{\vec{\theta}}}^{\otimes n}f,g)_{H^{\odot n}}
$$
will be a finite linear combination of the terms
$(a_i,\Pi_{\theta_j}h_k)_H$, hence the scalar product in $H^{\odot n}$
will be continuous with respect  
to  $\vec{\theta}\in [0,2\pi]^n$. Assume now that $(f_n)$ is a sequence of
finite linear combinations of pure vectors converging to $f$ in
$H^{\odot n}$. Then  
$$
\sup_{\vec{\theta}\in [0,2\pi]^n}
\left|\Bigl(f_k-f_l,\Pi_{\vec{\theta}}^{\otimes
    n}(h_1\odot\ldots\odot h_n)\Bigr)_{H^{\odot n}}\right|
\rightarrow 0\,,
$$
hence the limit is uniform, and this proves the continuity when $g$ is
a finite linear combination of pure vectors. Assume now that $g$ is
also a general symmetric tensor, then it can be approximated, as $f$,
by a sequence $(g_n)$ whose elements are the finite linear
combinations of pure vectors. Then we have again the following result 
\beaa
\lefteqn{\sup_{\vec{\theta}\in [0,1]^n}
\left|\Bigl(f,\Pi_{\vec{\theta}}^{\otimes n}(g_k-g_l)\Bigr)_{H^{\odot
      n}}\right|}\\ 
&&\leq \|f\|_{H^{\odot n}}\|g_k-g_l\|_{H^{\odot n}}\to 0\,,
\eeaa
which implies the uniform convergence with respect to $\vec{\theta}$. The
last claim is obvious when $f$ is a finite linear combination of pure
vectors. A  general $f$ can be approximated with such vectors, say
$(f_k,k\in \NN)$. Then, for any $x\in \reals^n$,  
$$
\int_{[0,2\pi]^n}e^{i(x,\vec{\theta})_{\R^n}}d(\Pi_{\vec{\theta}}^{\otimes
  n}f_k,f_k) 
=((R^{x_1}\otimes \ldots \otimes R^{x_n})f_k,f_k)_{H^{\odot n}}
$$
and this converges, as $k\to \infty$,  to the map 
$$
(x_1,\ldots,x_n)\to ((R^{x_1}\otimes \ldots \otimes
R^{x_n})f,f)_{H^{\odot n}}\,, 
$$
which is a continuous function on $\reals^n$ at $x=0$ by the spectral
representation of $R$. Then  the claim
follows from  the theorem of Paul L\'evy about the characterization of the weak
convergence of measures via the convergence of the characteristic functions.
\qed

We give now the main result of this section:
\begin{theorem}
\label{ergos-thm}
Let $R$ be a unitary operator on the Cameron-Martin space $H$ whose
resolution of identity is denoted by $(\Pi_\theta,\theta\in
[0,2\pi])$. Then the corresponding (measure preserving) transformation
$T$ is ergodic if and only if $\theta\to (\Pi_\theta h,k)_H$ is  
continuous on $[0,2\pi]$ for any $h,k\in H$. Moreover, if $T$ is ergodic,
it is also weak mixing.
\end{theorem}
\proof 
Let us first prove the necessity: assume that the resolution of
identity is not continuous. Then,  from Hahn-Banach theorem, there
exists an $h\in H$ and some $\tau\in (0,1)$ such that 
\beaa
(\Pi_{\tau+}h-\Pi_{\tau-}h,k)_H&=&\lim_{\eps\to
  0}(\Pi_{\tau+\eps}h-\Pi_{\tau-\eps}h,k)_H \\
&\neq& 0
\eeaa
for some $k\in H$. Let $z_\tau$ denote $\Pi_{\tau+}h-\Pi_{\tau-}h$.
Note that  we can represent  $z_\tau$ as 
$$
z_\tau=\int_{[0,2\pi]}\won_{\{\tau\}}(t)d\Pi_th\,.
$$
For any $k\in H$, the spectral representation of $R$ gives
\beaa
(Rz_\tau,k)_H&=&\lim_{\eps\to 0}\int_{[0,2\pi]}e^{i\theta}
d(\Pi_{\theta\wedge (\tau+\eps)}h-\Pi_{\theta\wedge(\tau-\eps)}h,k)\\
&=&e^{i\tau}(z_\tau,k)_H\,.
\eeaa
Therefore $z_\tau$ is an eigenfunction of $R$ with the corresponding
eigenvalue $e^{i\tau}$. Let $f(w)=|\delta z_\tau(w)|$. It is easy to
see that $f\circ T=f$ almost surely, hence $T$ can not be
ergodic and this contradiction proves the necessity. 
To prove the sufficiency, let $F$ 
be Wiener functional such that $F\circ T=F$ almost surely. Without
loss of generality we may assume that $F$ is bounded. From Lemma
\ref{chaos-lemma}, if we represent $F$ as $E[F]+\sum_n I_n(f_n)$, then
$I_n(f_n)\circ T=I_n(R^{\otimes n}(f_n))=I_n(f_n)$ for any $n\geq 1$.
Consequently
\beaa
0 &=& E\left[|I_n(f_n)-I_n(R^{\otimes n}(f_n))|^2\right]\\
&=& n!\left|f_n-R^{\otimes n}f_n\right|_{H^{\otimes n}}^2\\
&=& n!\int_{[0,2\pi]^n}\left|1-e^{i\sum_{k=1}^n \theta_k}\right|^2
 d((\Pi_{\theta_1}\otimes \ldots \Pi_{\theta_n})f_n,f_n)_{H^{\otimes n}}
\eeaa
this result implies that the positive measure 
$d((\Pi_{\theta_1}\otimes \ldots \Pi_{\theta_n})f_n,f_n)_{H^{\otimes
    n}}$ is concentrated on the set $\{\theta\in
[0,2\pi]^n:\exp i\sum_{1\leq k\leq n}\theta_k=1\}$, 
which is in  contradiction 
with the fact that it does not have atoms due to Lemma
\ref{conty-lemma}. The proof of weak mixing is similar with some
obvious modifications.

\qed

The mixing property of $T$ is straight forward.
\begin{theorem}
\label{mixing-thm}
The transformation $T$ defined by \req{4.1} is mixing if and only if 
$$
\lim_{n\to\infty} (R^nh,h)_H=0\,,
$$
for any $h\in H$.
\end{theorem}
\proof
By a density argument, $T$ is mixing if and only if 
$$
\lim_{n\to\infty}E\left[\rho(\delta h)\circ T^n\,\rho(\delta
  h)\right]=1
$$
for any $h\in H$, where 
$$
\rho(\delta h)=\exp\left\{\delta h-1/2|h|_H^2\right\}\,.
$$
We have 
\beaa
\rho(\delta h)\circ T^n\,\,\rho(\delta
h)&=&\exp\left\{\delta(R^nh+h)-|h|_H^2\right\}\\
&=&\rho(\delta(R^nh+h))\exp\left\{1/2|R^nh+h|_H^2-|h|_H^2\right\}\\
&=&\rho(\delta(R^nh+h))\exp(R^nh,h)_H\,.
\eeaa
Hence 
$$
\lim_{n\to\infty}E\left[\rho(\delta h)\circ T^n\,\rho(\delta
  h)\right]=1
$$
if and only if $(R^nh,h)_H\to 0$ as $n\to \infty$.
\qed

 We say that a sequence of random variables $(\eta_n,n\in \integers)$
 is ergodic or mixing if the shift transformation is ergodic or mixing 
 respectively. 
We have now the following corollary:
\begin{corollary}
The transformation $T$ is ergodic or mixing if and only if, for any
$h\in H$, the sequence $(\delta R^nh,\,n\in \integers)$ is ergodic or
mixing respectively.
\end{corollary}
\proof The necessity is evident, for the sufficiency it suffices to
remark that, Maruyama theorem implies the continuity of the spectral
measure associated to the sequence $(\delta R^nh,\,n\in \integers)$,
which is nothing but the measure $\theta \to d(\Pi_\theta h,h)_H$, whose 
continuity for any $h$ implies the ergodicity of $T$ by Theorem
\ref{ergos-thm}. For the mixing we proceed similarly.
\qed

\section{An example}

Let 
\beq{41}
dY_t = \gamma (t) d W_t; \qquad Y_0 = 0, t\in [0,1]
\eeq
where $W_\bfcdot$ is a standard $n$-dimensional Brownian motion and 
$\gamma(t), t\in [0,1]$ is a $n\times n$ unitary matrix, the
elements of $\gamma(t)$ will be assumed to be 
non-random and Lebesgue
measurable.  The ergodicity of $Y=T(w)$ will be discussed in this
section.
Let $e^{i\psi_j(t)}, 0<\psi_j \le 2\pi$, denote the eigenvalues of the
unitary matrix $\gamma$ and let $u(\cdot)$ denote the unit step function
$t$
$$
u(\alpha) = \left\{
\begin{array}{l l l}
1 & , & \alpha \ge 0\\
0 & , & \alpha < 0
\end{array}
\right. \,.
$$
Then we have
\begin{theorem}
A necessary and sufficient condition for the ergodicity of $T$ is the
continuity of
$\int_0^1 u(\theta-\psi_j(t))dt$ in $\theta \in [0, 2\pi]$ for all
$i=1, \cdots, n$.
Otherwise stated $T$ is ergodic iff the Lebesgue measure of
$C^j (\theta) = \{t: \psi_j(t,w)=\theta\}$ is zero for $\theta$ and all
$j$.
\end{theorem}

\proof
Assume first that for a.a.
$t\in [0, 1], \gamma$ possesses $n$-distant eigenvalues.  Also, assume
that $\psi_{j+1} > \psi_j$.
Since $\gamma(t)$ is unitary it has the representation
\beq{13}
\gamma(t) = A(t) \cdot \mathrm{diag.} e^{i\psi_j(t)} \cdot A^{-1} (t)
\eeq
where
$A(t) = [a_1 (t), \cdots a_n (t) ]$ is unique and
$\{a_j(t) j=1,2,\cdots, n\}$ are orthogonal $n$-vectors
$\gamma(t) \cdot a_i(t) = e^{i\psi_j (t)} a_j(t)$.
Let $h=\int_0^\bfcdot \dot{h}(s) ds$ where $\dot{h}$ takes values in 
$\reals^n$ and let $(\cdot , \cdot )$ denote the scalar product in
$\reals^n$ then
\beaa
\gamma(t)\dot{h} (t) & = &
\sum_j \Bigl(a_j(t), {h}' (t) \Bigr) e^{i\psi_j(t)} a_j(t) \\
& = & \int_0^{2\pi} e^{i\theta} d_\theta \sum_{j=1}^n u(\theta-\psi_j (t) \Bigr)
\Bigl(a_j (t), {h}'(t) \Bigr) \cdot a_j (t)
\eeaa
Hence
$$
Rh = \int_0^{2\pi} e^{i\theta} d_\theta \pi_\theta h
$$
where
$$
\Pi_\theta h = \int_0^\bfcdot \sum_j u\Bigl(\theta-\psi_j (t) \Bigr)
\Bigl(a_j(t), {h}'(t)\Bigr)
a_j (t) dt
$$
and
\beq{4.2}
|\Pi_\theta h|^2 = \int_0^{2\pi} \sum_j u\Bigl(\theta - \psi_j(t)\Bigr) 
\Bigl|
\Bigl(a_j(t), {h}'(t)\Bigr) \Bigl|^2 dt
\eeq
if, as $\eps\to 0$, 
$u(\theta+\eps-\psi_i(t)) \to u(\theta-\psi_i(t))$ for almost
(Lebesgue) all $t$ in $[0,2\pi]$ then by monotone convergence
$|\Pi_\theta h|^2$ is continuous in $\theta$.  Conversely,
if $|\Pi_\theta h|^2$ is discontinuous at $\theta=\theta_0$
$$
\lim_{\eps\to0}
\mathrm{Leb} \Bigl\{ t: u(\theta_0 + \eps - \psi_j (t) ) - u
(\theta_0-\eps, \psi_j(t)) \not=0 \Bigr\} > 0
\,.
$$ 
Hence $\int_0^{2\pi} u(\theta - \psi_j(t))dt$ is discontinuous at 
$\theta=\theta_0$. 
This proves the theorem for the case where there are $n$ distinct
eigenvalues.  If $\gamma$ possesses only $m<n$ distinct eigenvalues then $A(t)$
in \req{13} still holds with $\psi_j\le \psi_{j+1}$ but is no longer
unique.  This, however, can be overcome by constructing a measurable
selection which will provide a unique and measurable representation 
for $A(t)$.  The
rest of the proof remains unchanged.

\begin{corollary}
If $n$ is odd then $T$ is non ergodic.
\end{corollary}

\proof
If $n$ is odd then at least 
one of the eigenvalues of $\gamma(t)$ is either 1 
or -1. 
Now if $\psi_i (t)=\pi$ on a set of positive measure, 
then obviously \req{4.2} is discontinuous at
$\theta_0=\pi$.
For $\la_i=1$ on a set of positive $t$ measure, 
set $\psi_j = 2\pi$ (since $\Pi_\theta$ is, by definition
continuous as $\theta_2\searrow \theta_1$ and $\pi_0=0$)
and $\int_0^1 u(\theta-\psi_j(t)) dt$
must be discontinuous at $\theta_0 = 2\pi$.
\qed

\noindent
{\bf Acknowledgment:} The authors are grateful to 
L.~Decreusefond and
E.~Mayer-Wolf for
their  remarks.

\vspace{2cm}

{\footnotesize}
\begin{tabular}{ll}
A.S. \"Ust\"unel, & M. Zakai,\\
ENST,   D\'ept. R\'eseaux, \hspace{1.5cm }
& Department of Electrical Engineering,\\
46 Rue Barrault, & Technion---Israel Institute of Technology,  \\
75013 Paris&Haifa 32000,  \\
France&Israel\\
ustunel@enst.fr & zakai@ee.technion.ac.il
\end{tabular}

\end{document}